\documentclass[11pt]{article}

\usepackage{amsmath}
\usepackage{amssymb}
\usepackage{mathptmx}
\usepackage{amsfonts}
\usepackage[mathscr]{eucal}
\usepackage[dvips]{graphicx}
\usepackage{color}

\newcommand{\Yukawa}{}
\newcommand{\Dirac}{ \textrm{Dirac} }
\newcommand{\KG}{ \textrm{KG} }
\newcommand{\I}{\textrm{I}}

\newcommand{\fin}{\textrm{fin}}
\newcommand{\Dir}{ \textrm{Dir} }

\newcommand{\stateKG}{L^{2}( \mathbf{R}^{3}  )}
\newcommand{\stateDirac}{L^{2}( \mathbf{R}^{3} ; \mathbf{C}^{4} )}
\newcommand{\FYukawa}{\mathscr{F}}

\newcommand{\FDirac}{\mathscr{F}_{\textrm{Dirac}}}
\newcommand{\FKG}{\mathscr{F}_{\textrm{KG}}}

\newcommand{\FYukawaV}{\mathscr{F}_{V }}
\newcommand{\Yukawacore}{\ms{F}^{\textrm{fin}}_{\textrm{Dirac}} (\ms{D}(E))
 \hat{\otimes} \ms{F}^{\textrm{fin}}_{\textrm{KG}}  (\ms{D}(\omega))}
\newcommand{\YukawacoreVVprime}{\ms{F}^{\textrm{fin}}_{\textrm{Dirac}} (\ms{D}(E_{V}))
 \hat{\otimes} \ms{F}^{\textrm{fin}}_{\textrm{KG}}  (\ms{D}(\omega_{V'} ))}
\newcommand{\YukawacoreV}{\ms{F}^{\textrm{fin}}_{\textrm{Dirac}} (\ms{D}(E_{V}))
 \hat{\otimes} \ms{F}^{\textrm{fin}}_{\textrm{KG}}  (\ms{D}(\omega ))}

\newcommand{\HDirac}{H_{\textrm{Dirac}}}
\newcommand{\HKG}{H_{\textrm{KG}}}

\newcommand{\Hzero}{H_{0}}

\newcommand{\HDiracV}{H_{\textrm{Dirac} ,V}}

\newcommand{\mbf}[1]{\ensuremath{\mathbf{#1}}}
\newcommand{\ms}[1]{\ensuremath{\mathscr{#1}}}
\newcommand{\sqz}[1]{\ensuremath{d\Gamma({#1}) }}

\newcommand{\tens}{\otimes}

\newcommand{\nstens}{\otimes^{n}_{s}}
\newcommand{\natens}{\otimes^{n}_{a}}

\newcommand{\eltwo}{L^{2}(\mathbf{R}^{3})}
\newcommand{\Rthree}{\mathbf{R}^{3} }
\newcommand{\dx}{d \mathbf{x} }

\newcommand{\dk}{d \mathbf{k} }

\newcommand{\restr}{\upharpoonright}

\newcommand{\psivarxL}{ \overline{\psi_{L} (\mbf{x}) } }
\newcommand{\psivarxLV}{ \overline{\psi_{L,V} (\mbf{x}) } }

\newcommand{\cutoff}[2]{\ensuremath{
\chi_{\textrm{#1}} (\mathbf{#2} )}}

\newtheorem{theorem}{Theorem}[section]
\newtheorem{proposition}[theorem]{Proposition}
\newtheorem{lemma}[theorem]{Lemma}

\begin{document}

\title{ Ground States of the  Yukawa Models with Cutoffs
}
\author{Toshimitsu TAKAESU}
\date{ }
\maketitle
\begin{center}
\textit{Faculty of Mathematics, Kyushu University, \\
 Fukuoka, 812-8581, Japan }
\end{center}

\begin{quote}
\textbf{Abstract.}
       Ground states of the
           so called Yukawa model is considered.  The Yukawa model describes a Dirac field interacting with a Klein-Gordon field. By introducing both ultraviolet cutoffs  and  spatial cutoffs, the total Hamiltonian is defined as a self-adjoint operator on a boson-fermion Fock space. It is shown that the total Hamiltonian has  a positive spectral gap
               for \textit{all} values of coupling constants.  In particular the existence of ground states is proven.                 \\
$\; $ \\
{\small
key words : Spectral analysis, Fock spaces, Quantum field theory}.
\end{quote}

\section{Introduction}
$\; \; $
In this paper we investigate the existence of ground states of the Yukawa model which
describes a Dirac field  interactiong with a Klein-Gordon field.
 Both Dirac field and
 Klein-Gordon field are massive,
 and ultraviolet cutoffs are imposed on both of them.
 The total Hamiltonian of the Yukawa model
 is the sum of the free Hamiltonian
    and the interaction Hamiltonian:
\begin{equation}
H \; \; = \;  \;  \HDirac\otimes I
+ I \otimes  \HKG    \; + \kappa \,  H' \end{equation}
on
$ \FYukawa \; =   \FDirac \tens \FKG$, where $\kappa \; > \; 0$ is  a coupling constant.
The free Hamiltonians $\HDirac$ and $\HKG$ are given by formally
\begin{align*}
&\HDirac \; = \sum_{s=\pm 1/2}\int_{\Rthree}
 \sqrt{M^2 + \mbf{p}^2} \left( \frac{}{}
  b_{s}^{\ast} (\mbf{p}) b_{s} (\mbf{p}) \; + \;  d_{s}^{\ast} (\mbf{p}  ) d_{s} (\mbf{p}) \right) d \mbf{p} ,  \qquad \qquad M > 0 , \\
  & \HKG \; = \int_{\Rthree}
 \sqrt{m^2 + \mbf{k}^2} a^{\ast}(\mbf{k})  a(\mbf{k}) \dk , \qquad \qquad \qquad m>0 .
\end{align*}
In the subsequent section, we give the rigorous definition of   $\HDirac$ and $\HKG$.
The interaction Hamiltonian
$H'$
is defined by
\begin{equation}
H' \; \; = \; \; \int_{\Rthree} \cutoff{I}{x}
 \overline{\psi_{{}_{\chi_{{}_{\Dir}}}} (\mbf{x})} \psi_{{}_{\chi_{{}_{\Dir}}}} (\mbf{x}) \tens
  \phi_{{}_{\chi_{{}_{\KG}}}} (\mbf{x} ) \dx , \notag
\end{equation}
where   $\psi_{\chi_{{}_{\Dir}}} (\mbf{x})$ is
a Dirac field
with   an  ultraviolet cutoff
$ \chi_{{}_{\Dir}}$, and  $ \phi_{\chi_{{}_{\KG}}} (\mbf{x} ) $  a Klein-Gordon field with  an ultraviolet cutoff $\chi_{{}_{\KG}}$.
We furthermore introduce a spatial cutoff  $ \cutoff{I}{x} $ in $H'$
to  define $H$ as  a  self-adjoint operator.
 Since
 the interaction $H'$ is
 relatively bounded
 with respect to
 $\HDirac \otimes 1+ 1\otimes \HKG$ by virtue of
 cutoffs,
 $H$ is self-adjoint and bounded from below by  the Kato-Rellich theorem.
 We say that a self-adjoint  operator $X $ bounded from below has a  ground state,
 if  the bottom of its spectrum is an eigenvalue, and
 the difference between the bottom of the spectrum and that of the essential spectrum is called spectral gap.
In this paper we show that  $H$ has
 a positive spectral gap
   for all values of coupling constants. In particular the existence of ground states follows from this.
 $\quad$ \\
$\quad$ In the last decade,  a system of
quantum particles governed by a Schr\"odinger operator
 interacting with a massless bose field are successfully investigated.
In particular
the existence of ground states
of some massless models in non-relativistic
QED  is proven in \cite{Ge00, GLL01} for all values of coupling
constants.
 It is also shown  in \cite{AGG07, BG09, BDG04, Ta09} that ground states of
 a massless model in QED
      exist but for sufficiently small values of
  coupling constants.
 Since quantized radiation fields in QED and in nonrelativistic QED are massless,
 the spectral gap of the free Hamiltonians is zero.
   Then all the results mentioned above are not trivial.
     For other topics on the system of fields interacting fields, refer to \cite{AGG07, BG09}.  On the analysis of a   field equation of the Yukawa model, called the Dirac-Klein-Gordon equation, see \cite{NB00,SM07,RBEW}.

$\quad $ \\
$\quad $ Now let  us consider
the existence of ground states of the Yukawa model
$H$.
 Since $\HDirac$ and $\HKG$ are massive,
the spectral gap of
$\HDirac \otimes 1+1\otimes  \HKG$ is positive.
Then the regular perturbation  theory \cite{RS}
 says that
   $H$ also has ground states for sufficiently small values of coupling constants.
 It is not obvious, however,  whether $H$ also has ground states  for all values of coupling constants.
Moreover unfortunately we can not directly apply methods developed in  \cite{Ge00,GLL01} to show the existence of ground states of $H$. Outline of our strategy is as follows. To prove the existence of ground states of $H$,
we  use
a momentum lattice approximation \cite{GJ70, AH97}.
Then $H$ can be approximated with
   some lattice parameters $V$ and $L$ as
$$H_{L,V}=\HDiracV \otimes 1 +1\otimes  \HKG+\kappa H_{L,V}'.$$
 It is shown that
 $\HDiracV$ has a  compact resolvent. Then from a standard argument as in
   \cite{Ge00,AH97}, it follows that
    $H_{L,V}$ has a positive  spectral gap which is uniform with respect to $V$ and $L$ by positive masses $m$ and $M$.
 Since $H_{L,V}$ converges to $H$ in the uniform resolvent sense as $V\to\infty$ and $L\to\infty$,
       we can see that $H$ also has a positive spectral gap.
In this paper
integrable condition
$\int_{\Rthree} |\mbf{x} | \, |  \cutoff{I}{x} | d \mbf{x}  < \infty $ is supposed.  This assumption corresponds
to
the spatial localization discussed in \cite{BFS99,GLL01}.

$\quad $ \\
$\quad$ \\
$\quad$  This paper is organized as follows.
In Section 2, we introduce Dirac fields and Klein-Gordon fields with ultraviolet cutoffs. Then we  define the Yukawa Hamiltonian with spatial cutoffs on a boson-fermion Fock space, and state  a main result.
In Section 3, we give the proof of the  main theorem.

\section{Definitions and Main Results}
\subsection{Dirac Fields and Klein-Gordon Fields}
We first consider Dirac fields.
The state space defined by
$ \ms{F}_{\Dirac} = \oplus_{n=0}^{\infty}
 (  \natens \stateDirac ) $, where $\natens \stateDirac $ denotes the $n$-fold anti-symmetric tensor product
of $ \stateDirac $ with
$\; \tens_{a}^{0} L^{2}(\Rthree ; \mbf{C}^{4} ) ~:= \mbf{C}$.
Let $\ms{D} $ be the subset of $\stateDirac $.
 We define  the finite particle subspace
 $\ms{F}_{\Dirac}^{\; \fin} (\ms{D})$ on $\ms{D}$ by the set  of $\Psi = \{\Psi^{(n)} \}_{n=0}^{\infty} \; $ satisfying  that
   $\Psi^{(n)}  \in \natens \ms{D}  $ and
    $\Psi^{(n')}  = 0  $ for all $ n' > N$ with  some $N \geq 0$.
Let $B(\xi )$, $ \; \xi ={}^{t} (\xi_{1}, \cdots , \xi_{4} )  \in \stateDirac$, and
$\; B^{\ast}(\eta)$, $ \; \eta ={}^{t} (\eta_{1}, \cdots , \eta_{4} )  \in\stateDirac$, be the
 annihilation operator and the creation operator on  $\FDirac$, respectively.
For $f \in \eltwo$ let us set
\begin{align*}
&b^{\ast}_{1/2}(f) = B^{\ast}({}^{t} (f,0,0,0) ), \quad \quad
b^{\ast}_{-1/2}(f) = B^{\ast}({}^{t}(0,f,0,0) ), \\
&d^{\ast}_{1/2}(f) = B^{\ast}({}^{t}(0,0,f,0) ), \quad \quad
d^{\ast}_{-1/2}(f) = B^{\ast}({}^{t}(0,0,0,f) ) .
\end{align*}
Then they satisfy canonical anti-commutation relations :
 \begin{align*}
&\{ b_{s} (f) , b_{\tau}^{\ast}(g) \} =
\{ d_{s} (f) , d_{\tau}^{\ast}(g) \} =
\delta_{s, \tau} (f, g)_{\eltwo}   ,  \\
&\{ b_{s} (f) , b_{\tau} (g) \} =
\{ d_{s} (f) , d_{\tau} (g) \} = \{ b_{s} (f) , d_{\tau} (g) \} =
\{ b_{s} (f) , d_{\tau}^{\ast}(g) \} = 0  .
\end{align*}
It is known that $b_{s} (\xi )$ and $d_{s} (\xi )$ are bounded with
\begin{equation}
 \| b_{s} (\xi )  \|  = \| d_{s} (\xi )  \| = \| \xi  \| .  \label{boundbd}
 \end{equation}
The one particle energy of Dirac field with momentum $\mbf{p} \in \Rthree $ is given by
 $ \;
E (\mbf{p})  = \sqrt{\mbf{p}^2  \; + \; M^2 }$,
where $M>0 $ denotes the mass of
an electron.
Let
  \begin{equation}
f_{s}^{l} (\mbf{p}) =
\frac{\chi_{\Dir}(\mbf{p} ) u_{s}^{l} (\mbf{p}) }{ \sqrt{ (2 \pi ) ^3  E (\mbf{p})}}     , \qquad
g_{s}^{l} (\mbf{p})=
\frac{\chi_{\Dir} (\mbf{p})    v_{s}^{l} (-\mbf{p})}
{\sqrt{ (2 \pi ) ^3  E (\mbf{p} )}} , \qquad  \qquad s= \pm 1/2 , \qquad  l= 1, \cdots  ,4 , \notag
\end{equation}
where  $\chi_{\Dir}$ is an ultraviolet cutoff, and  $ u_{s}(\mbf{p} ) =
(u_{s}^{l} (\mbf{p}) )_{l=1}^{4} \; $ and  $v_{s}(\mbf{p} ) = (v_{s}^{l} (\mbf{p}) )_{l=1}^{4}  \; $ denote spinors with the positive and negative energy part of $\mbf{\alpha} \mbf{\cdot}
  \mbf{p}   + \beta M $  with spin $s =  \pm 1/2 $, respectively. Here $\alpha^{j} $, $j=1,2,3$, and  $\beta$ are the $4\times 4$ matrix satisfying the canonical anti-commutation relation $
 \{ \alpha^{j} , \alpha^{l} \} = 2 \delta_{j,l} , \;
\{ \alpha_{j} , \beta \} = 0 ,   \; \beta^{2}=I  $.
The Dirac field
$\; \psi (\mbf{x}) \, = \, {}^{t}( \psi_{1}(\mbf{x}), \cdots , \psi_{4}(\mbf{x}) ) $
   is defined by
\[
\psi_{l}(\mbf{x}) = \sum_{s=\pm 1/2}( b_{s} (f_{s,\mbf{x}}^{l} ) +   d^{\ast}_{s} (g_{s, \mbf{x}}^{l} )) ,
 \qquad \qquad   l= 1, \cdots  ,4 ,
\]
where $\;   f_{s,\mbf{x}}^{l} (\mbf{p}) = f_{s}^{l} (\mbf{p})
e^{-i \mbf{p} \cdot \mbf{x}} \; $ and $ \;  g_{s,\mbf{x}}^{l} (\mbf{p}) = g_{s}^{l} (\mbf{p})
e^{-i \mbf{p} \cdot \mbf{x}} $.
We introduce   the following assumption.
\begin{quote}
\textbf{(A.1) (Ultraviolet cutoff for  Dirac fields) }
$\; $ $\chi_{\Dirac}$ satisfies that
\[
\int_{\Rthree} \frac{ | \chi_{\Dir} (\mbf{p})
  u_{s}^{l} (\mbf{p})|^2 }{  E_{M}(\mbf{p}) }  d \mbf{p}  \; <  \; \infty  , \qquad
\int_{\Rthree}  \frac{ | \chi_{\Dir} (\mbf{p})
  v_{s}^{l} (- \mbf{p})|^2 }{  E_{M}(\mbf{p}) } d \mbf{p}  \; <  \; \infty .
\]
\end{quote}
$\quad$ \\

We secondly define Klein-Gordon fields. The state space is defined by $ \FKG \;   = \; \oplus_{n=0}^{\infty}
 (  \nstens \stateKG ) $,
where $\nstens \stateKG$ denotes the $n$-fold symmetric tenser
 product of $\stateKG \; $ with $\tens_{s}^{0} \stateKG := \mbf{C}$.
  In a similar way to the case of  Dirac fields,
   we define the finite particle subspace
 $\ms{F}_{\KG}^{\; \fin} (\ms{M})$  on $\ms{M}
 \subset \stateKG$ but
   anti-symmetric tensor products is replaced by
    symmetric tensor products.
  Let $a( \xi ), \; \xi  \in \stateKG $, and $a^{\ast}(\eta ), \; \eta \in \stateKG $, be the annihilation operator and the creation operator on $\FKG$, respectively.
  Then they satisfy canonical commutation relations on  $\ms{F}_{\KG}^{\fin}(  \stateKG )$ :
\[
    [ \, a( \xi ), \, a^{\ast} (\eta )  ] = (\xi , \eta  ), \qquad
    [ \, a  (\xi ), \, a (\eta  )  ] = [ a^{\ast}(\xi  ), \,  a^{\ast}(\eta ) ] =0 .
\]
Let $S $ be a  self-adjoint operator on $\stateKG$.
The  second quantization of $S$ is defined by
\[
 \sqz{S}_{\restr {}_{\ms{F}_{\KG}}}  \; = \bigoplus_{n=0}^{\infty} \left(
\sum_{j=1}^{n} ( I \tens \cdots I \tens \underbrace{S}_{jth} \tens I   \cdots  \tens I )\right)_{\restr {}_{\ms{F}_{\KG}}} .
\]
Similarly, we can define the second quantization 
 $\sqz{A}_{\restr {}_{\ms{F}_{\Dirac}}} $ of the Dirac field for a  operator $A$ on $\stateDirac$.
For $\eta \in \ms{D} (S^{-1/2})$, $a(\eta) $ and $a^{\ast} (\eta)$ are  relatively bounded with respect to $\sqz{S}_{\restr {}_{\FKG}} $ with
\begin{align}
& \| a (\eta ) \Psi  \| \leq \| S^{-1/2} \eta  \|  \, \|  \sqz{S}^{1/2}_{\restr {}_{\FKG}} \Psi \| ,
    \qquad \qquad  \qquad  \quad  \Psi  \in  \ms{D} ( \sqz{S}_{\restr {}_{\FKG}}^{1/2} ),    \label{bounda}   \\
& \| a^{\ast}(\eta ) \Psi  \| \leq \| S^{-1/2} \eta \|
\|   \sqz{S}^{1/2}_{\restr {}_{\FKG}} \Psi \| + \| \eta  \| \| \Psi \| ,  \qquad   \Psi  \in  \ms{D} (
 \sqz{S}_{\restr {}_{\FKG}}^{1/2} ) .
\label{boundad}
\end{align}
The one particle energy of Klein-Gordon field with momentum $\mbf{k} \in \Rthree$ is given  by
$ \omega (\mbf{k})  = \sqrt{\mbf{k}^2  \; + \; m^2 } $, $ m> 0 $.
Let us define the field operator $\phi (\mbf{x}) $ by
\[
 \phi (\mbf{x}  ) \; = \; \frac{1}{\sqrt{2}}
 \left( \frac{}{} a (h_{\mbf{x}}) \; + \;  a^{\ast} (h_{\mbf{x}}) \right),
\]
where $h_{\mbf{x}}  (\mbf{k})  = h (\mbf{k}) e^{i\mbf{k} \cdot \mbf{x}} \; $ with
$ \; h (\mbf{k}) \; = \; \frac{\cutoff{KG}{k} }{\sqrt{(2 \pi )^3 \omega  (\mbf{k}) }}$,
and $\chi_{\KG}$ is an ultraviolet cutoff function. We assume  the following condition :
 \begin{quote}
\textbf{(A.3)}  \textbf{(Ultraviolet cutoffs for Klein-Gordon fields)}
 $\quad$ $\chi_{\KG}$ satisfies that
\[
\int_{\Rthree}
 \frac{ | \cutoff{KG}{k} |^2 }{ \omega  (\mbf{k}) }
 < \infty ,  \qquad \qquad \int_{\Rthree}
   \frac{ | \cutoff{KG}{k} |^2 }{ \omega  (\mbf{k})^2 }
 < \infty .
 \]
\end{quote}

\subsection{Total Hamiltonian and Main Theorem}
The state space of the interaction system between  Dirac fields and  Klein Gordon fields
 is given by
 \[
 \FYukawa \; =   \: \FDirac \tens \FKG ,
 \]
and the free Hamiltonian by
\[
\Hzero \; = \HDirac  \tens I \; + \; I \tens \HKG ,
\]
where
$  \HDirac \; = \;    \sqz{E}_{\restr \FDirac}  $    and
$ \HKG \; = \; \sqz{\omega}_{\restr \, \FKG}$.
To define the interaction, we  introduce a spatial  cutoff
 satisfying the  following condition :
 \begin{quote}
\textbf{(A.3)}  \textbf{(Spatial cutoffs)}
 $\;$ $\chi_{\I}$ satisfies that
$
\int_{\Rthree} | \chi_{\I} (\mbf{x}) | \dx < \infty$  .
\end{quote}
Now let us define the linear functional $   \FYukawa  \, \times  \, \left( \ms{F}_{\Dirac}^{\fin}
(\ms{D} (E) ) \hat{\tens}   \ms{F}_{\KG}^{\fin}
(\ms{D} (\omega ) ) \frac{}{} \right)  \to \mbf{C} $, where $\hat{\tens}$ denotes the algebraic tensor product, by
\begin{equation}
\ell_{\I} (\Phi , \Psi  ) \; = \;
\int_{\Rthree} \cutoff{I}{x} \left( \Phi , \,
 \overline{\psi (\mbf{x}) } \psi (\mbf{x} ) \tens \phi (\mbf{x})
 \Psi  \right) \dx ,
 \end{equation}
where
$ \;  \overline{\psi (\mbf{x}) } \, = \, \psi^{\ast}(\mbf{x} ) \gamma^0$
with $\gamma^0=\beta$.
By (\ref{boundbd}) we have
\begin{equation}
\| \psi_{l}  (\mbf{x} ) \| \; \leq  \;  M_{\Dir}^{l} , \label{bound_psil}
\end{equation}
where $  M_{\Dir}^{l} = \sum_{s=\pm 1/2 } (   \| f_{s}^{l} \| +   \| g_{s}^{l} \|  ) $. We also see that
by (\ref{bounda}) and (\ref{boundad}),
\begin{equation}
\| \phi_{l}  (\mbf{x} )  \Psi \| \;
\leq \sqrt{2} M_{\KG}^{1} \| \HKG^{1/2} \Psi \|
 \; + \; \frac{1}{\sqrt{2}} M_{\KG}^{0} \| \Psi \| , \label{bound_phi}
\end{equation}
where $ M_{\KG}^{j} =  \|  \frac{h}{\sqrt{\omega^j } } \| $, $\; j \; \in  \{ 0 \} \cup  \mbf{N}$.
By (\ref{bound_psil}) and (\ref{bound_phi}), we have
\begin{equation}
| \ell_{\I} ( \Phi , \Psi ) | \; \leq \; \left(  L_{\I}  \| (I \otimes H_{\KG}^{1/2} ) \Psi \| + R_{\I} \| \Psi \| \right) \| \Phi \|  ,   \label{ell_I_bound}
\end{equation}
where
 $ L_{\I} = \sqrt{2} \| \chi_{\I} \|_{L^{1}} \sum_{l, l' }
| \gamma_{l, \, l'}^{0} | \,
\, M^{\, l}_{\Dir} \, M^{\, l'}_{\Dir} \, M^{\,1}_{\KG} $, and $\;  R_{\I} = \frac{1}{\sqrt{2}} \| \chi_{\I} \|_{L^{1}} \sum_{l, l' }
| \gamma_{l, \, l'}^{0} | \,
\, M^{\, l}_{\Dir} M^{\, l'}_{\Dir} \, M^{\,0}_{\KG} $.
By the Riesz representation theorem,  we can define the symmetric operator  $H' \, : \, \ms{F}_{\Yukawa} \; \to  \ms{F}_{\Yukawa} $ such  that
\begin{equation}
   ( \Phi , H' \Psi   ) \; \; =  \; \; \ell_{\I}( \Phi, \Psi  )     ,  \label{def_interactI}
\end{equation}
 and
\begin{equation}
\|  H'  \Psi  \|  \; \leq  \;   L_{\I}  \| (I \otimes H_{\KG}^{1/2} ) \Psi \| + R_{\I} \| \Psi \| .
\label{H_I_bound}
\end{equation}
We see that $ H' $ is formally denoted by
\[
 H'  \; = \;  \int_{\Rthree} \cutoff{I}{x}
 \overline{\psi (\mbf{x}) } \psi (\mbf{x} ) \tens \phi (\mbf{x}) \dx .
\]
The total Hamiltonian of the Yukawa model is then defined by
\begin{equation}
\qquad H \; \;  = \; \; \Hzero \; +  \kappa \, H' , \qquad \qquad  \kappa \in \mbf{R}.
\end{equation}
Let us consider the self-sdjointness of $H$.
For  $\epsilon \; > \; 0$, there exists  $C_{\epsilon} \geq 0 $ such that
 for all $\Psi \in \ms{D}(H_{\KG})$,
\begin{equation}
\| H_{\KG}^{1/2} \Psi \| \leq
\epsilon \| H_{\KG} \Psi \| + c_{\epsilon} \| \Psi \| .
\label{halfepsilon}
\end{equation}
  Then by (\ref{halfepsilon}) and (\ref{H_I_bound}), we see that for $\Psi \in \ms{D} (H_{0})$,
\begin{equation}
\| H' \Psi \| \leq \epsilon L_{\I} \| H_{0} \Psi \|  +
(c_{\epsilon} L_{\I} +R_{\I})
\| \Psi \| . \label{8/30.1}
\end{equation}
Let us  take  sufficiently small $\epsilon >0 $ such as $ \epsilon L_{\I} \; < \; 1$ in (\ref{8/30.1}). Then  by the Kato-Rellich theorem,
$ H $ is self-adjoint on  $\ms{D} (H_{0})$ and essentially self-adjoint on any core of $H_{0}$.
In particular, $H $ is essentially self-adjoint on
\begin{equation}
\ms{D}_{0} \; = \; \Yukawacore .  \label{Yukawa_core}
\end{equation}
The Kato-Rellich theorem also shows that $H$ is bounded from below i.e.
 $\inf \sigma (H) \; > - \infty$.

$\quad $ \\
$\quad$
Let $X$ be  self-adjoint and
 bounded from below.
Let  us denote the infimum of the spectrum of $X$ by $E_{0} (X) \; = \; \inf \sigma (X)  $.
We say that  $X $ has a ground state if
$\; E_{0} (X) $ is an eigenvallue of $X$. \\
 Let
\begin{equation}
\nu \; \; =  \; \;  \min  \;  \{  m, \; M \} .
\end{equation}
Then it is known that  the spectrum of $\Hzero$ is $\sigma (H_{0}) = \{ 0 \} \cup [ \nu , \infty )$.
 To prove the existence of the ground states of $H$,  we introduce the additional condition on the spatial cutoff.
\begin{quote}
\textbf{(A.4) (Spatial localization) } 
 $\;$ $\chi_{\I}$ satisfies that
$
\int_{\Rthree} |\mbf{x} | \, | \chi_{\I} (\mbf{x}) | \dx < \infty $.
\end{quote}
Now we are in the position to state the main theorem.
\begin{theorem}  \label{MainTheorem} $\, $ \\
Assume (\textbf{A.1})-(\textbf{A.4}).
Then $[E_0(H),E_0(H)+\nu)\cap \sigma(H)$ is purely discrete  for all values of  coupling constants.
 In particular $H$ has  ground states for all values of  coupling constants.
\end{theorem}

\section{Proof of Main Theorem}
Let us introduce some notations. Let $\Gamma_{V}$ be the set of lattice points
\[
\Gamma_{V} \;
=\; \{ \mbf{q} = (q_{1},q_{2},q_{3}) \;  | \;  q_{j} = \frac{2\pi}{V} n_{j}  ,\; \;  n_{j} \in \mbf{Z} , \;  j=1,2,3  \} .
\]
For  each lattice point $\mbf{q} \in \Gamma_{V} $, set
$
C ( \mbf{q} ,V ) \;  = \;  [ q_{1}-\frac{\pi}{V}, q_{1}+\frac{\pi}{V} )  \times
 [ q_{2}-\frac{\pi}{V}, q_{2}+\frac{\pi}{V} ) \times
  [ q_{3}-\frac{\pi}{V}, q_{3}+\frac{\pi}{V} ) \; \;
\subset \; \Rthree  $
and   $ I_{L} =  [-L , L ]  \times [-L , L ]  \times [-L , L ]  \; \subset \Rthree $.
For $\xi \in \eltwo $,
we define the approximated  functions $ \xi_{L}$ and $\xi_{L,V}$ by
\begin{align*}
&  \xi_{L} (\mbf{k} ) =\xi (\mbf{k} )  \chi_{I_{L}}  (\mbf{k} ),  \\
& \xi_{L,V} (\mbf{k} ) = \sum_{\mbf{q} \in \Gamma_{V}}   \xi
 (\mbf{q} ) \chi_{C( \mbf{q}, V)  \cap I_{L} }(\mbf{k} )   ,
\end{align*}
where  $\chi_{J} (\mbf{k}) $ denotes the characteristic function on $J \subset \Rthree$.
By considering  the  map  $ \stateKG \ni \xi $ $ = \sum_{\mbf{q} }  \xi
 (\mbf{q} ) \chi_{C( \mbf{q}, V)}  $ $ \mapsto   (\xi (\mbf{q}))_{\mbf{q} \in \Gamma_{V}}  \in \ell^{2}( \Gamma_{V})$, we can identify  $ \ell^2 (\Gamma_{V} ) $ as a closed subspace of
 $\stateKG $.
 Let us set
\[
\ms{F}_{\Yukawa  V}  \; = \;   \ms{F}_{\Dirac , V } \tens \FKG ,
\]
where $
 \ms{F}_{\Dirac , V}  \; = \; \oplus_{n=0}^{\infty}
 (  \natens  \ell^2 (\Gamma_{V}   ;  \mbf{C}^4 ) )   $.
Let us
define $H_{0,V}$ on
 $\ms F_{\Yukawa}$ by
\[
H_{ 0, V} \;  =  \;
\HDiracV \tens I + I \tens \HKG ,
\]
where $ \HDiracV = \sqz{E_V }_{\restr \ms{F}_{\Dirac }}$ with $ E_{V} (\mbf{p} ) = \sum\limits_{\mbf{q} \in \Gamma_{V}} E(\mbf{q} )
\chi_{C  ( \mbf{q} , V) }(\mbf{p} )$.
Approximated interaction Hamiltonians are also defined by
\begin{align*}
 & H_{L,V}'  \int_{\Rthree} \chi_{\I}(\mbf{x} )
 \left( \psivarxLV \psi_{L,V} (\mbf{x} ) \otimes  \phi (\mbf{x} ) \right) \dx , \\
 & H_{L}' =  \int_{\Rthree} \chi_{\I}(\mbf{x} )
 \left( \psivarxL \psi_{L} (\mbf{x} ) \otimes  \phi (\mbf{x} ) \right) \dx ,
\end{align*}
where $\psi_{L} (\mbf{x}  )  = ( \psi_{L}^l (\mbf{x}  ) )_{l=1}^4 \; $  and
$\;  \psi_{L,V } (\mbf{x}  )  = ( \psi_{L,V}^l (\mbf{x}  ) )_{l=1}^4 $ with
  $ \; \psi_{L}^l (\mbf{x}) = \sum\limits_{s=\pm 1/2}
\{ b_{s} (( f_{s,\mbf{x}}^{l} )_L )  $  $ +   d^{\ast}_{s} (( g_{s, \mbf{x}}^{l}  )_{L}  )  \} \;   $ and
$\;  \psi_{L,V}^l (\mbf{x}) = \sum\limits_{s=\pm 1/2}
\{ \frac{}{} b_{s} ( ( f_{s,\mbf{x}}^{l}  )_{L,V} ) $ $ +   d^{\ast}_{s} ((g_{s, \mbf{x}}^{l}  )_{L,V}  )  \}  $.
 Let
\begin{align}
&H_{L,V} \;  = \; \; H_{0,V} \; + \; \kappa H_{L, V}' ,   \\
&H_{L} \;\; \; \;  = \; \; H_{0} \; + \; \kappa H_{L}'  .
\end{align}
In a similar way to the case of $H$, we can prove that $H_{L} $ and $H_{L,V} $ are  essentially self-adjoint on $\ms{D}_{0}$ and $\ms{D}_{0,V} = \YukawacoreV$, respectively.
\begin{lemma} \label{11/8.1}
Assume \textbf{(A.1)}-\textbf{(A.3)}.
Then  $\; H_{L,V}$ is reduced to
$\ms{F}_{\Yukawa  V} $.
\end{lemma}
\textbf{(Proof)}
Let us denote $p_{V}$ the orthogonal projections from $\eltwo $
  to $ \ell^{2}(\Gamma_{V})  $. Then $\Gamma (p_{V}) = \oplus_{n=0}^{\infty} (\otimes^n p_{V})$ is the projection
   from $\FDirac$ to $ \ms{F}_{\Dirac , V}$.
 Let $\Psi  \in \ms{D}_{0,V} $.
 Then it is easy to see that
$( \Gamma (p_{V}) \tens I  )   H_{0,V} \Psi =  H_{0,V}
  ( \Gamma (p_{V}) \tens I  )  \Psi \; $.
By using $p_{V} \chi_{C(\mbf{q},V)} = \chi_{C(\mbf{q},V)} $,  we  see that for all $\Phi \in \ms{F}_{\Yukawa}$,
     \[
 \int \chi_{\I}(\mbf{x} )
( \Phi, ( \Gamma (p_{V}) \overline{\psi_{L,V} (\mbf{x})}  \psi_{L,V} (\mbf{x} ) )
 \otimes   \phi (\mbf{x} )  \Psi ) \dx
=  \int\chi_{\I}(\mbf{x} )
( \Phi,   (\overline{\psi_{L,V} (\mbf{x})}  \psi_{L,V} (\mbf{x} )
\otimes  \phi   (\mbf{x} )  )(\Gamma (p_{V}) \otimes I ) \Psi ) \dx.
  \]
 Hence $ ( \Gamma (p_{V})  \otimes I  )   H_{L,V}' \Psi  =
 H_{L,V}' (  \Gamma (p_{V})  \otimes I )  \Psi $.  Thus
$   ( \Gamma (p_{V})  \otimes I  )   H_{L,V} \Psi =  H_{L,V}
  ( \Gamma (p_{V})  \otimes I  )   \Psi \; $ follows for all  $\Psi  \in \ms{D}_{0,V} $. Since
 $ \ms{D}_{0,V} $ is a core of $H_{L,V}$,
the lemma follows.  $\blacksquare $    

\begin{proposition}   \label{11/8.2}
Assume \textbf{(A.1)}-\textbf{(A.4)}. Then
$H_{L,V \, \restr_{\ms{F}_{\Yukawa  V}}}$ \text{ has  purely discrete spectrum  in }
$[ E_{0}( H_{L,V}  ),  E_{0}( H_{L,V}   )+ \nu  )$.
\end{proposition}

$\quad$ \\
To prove Proposition \ref{11/8.2}, we also take the lattice approximation of Klein-Gordon fields.
Let us set
\[
\ms{F}_{\Yukawa  V,V'}  \; = \;   \ms{F}_{\Dirac , V } \tens \ms{F}_{\KG ,V'} ,
\]
where $
 \ms{F}_{\KG , V'}  \; = \; \oplus_{n=0}^{\infty}
 (  \nstens  \ell^2 (\Gamma_{V'}   ))   $.
Set
\[
H_{ 0, V,V' } \;  =  \;
\HDiracV \tens I + I \tens H_{\KG , V'} ,
\]
where $ H_{\KG , V'}  =  \sqz{\omega_{V'} }_{\restr \ms{F}_{\KG}}$ with $ \omega_{V'} (\mbf{k} ) = \sum\limits_{\mbf{q} \in \Gamma_{V'}} \omega (\mbf{q} )
\chi_{C  ( \mbf{q} , V') }(\mbf{k} )$. Let
\begin{align*}
& H_{L,V,L'}' =  \int_{\Rthree} \chi_{\I}(\mbf{x} )
 \left( \psivarxLV \psi_{L,V} (\mbf{x} ) \otimes  \phi_{L'} (\mbf{x} ) \right) \dx , \\
 & H_{L,V,L' ,V'}'  \int_{\Rthree} \chi_{\I}(\mbf{x} )
 \left( \psivarxLV \psi_{L,V} (\mbf{x} ) \otimes  \phi_{L' , V'} (\mbf{x} ) \right) \dx ,
\end{align*}
where
$ \phi_{L'} (\mbf{x} ) = \frac{1}{\sqrt{2}}
 \left\{ \frac{}{} a ((h_{\mbf{x} })_{L'} ) +
 a^{\ast}((h_{ \mbf{x} })_{L'} )\right\} \; $  and
   $  \;     \phi_{L',V'} (\mbf{x} ) = \frac{1}{\sqrt{2}}
 \left\{ a ((h_{\mbf{x} })_{L',V'} ) +
 a^{\ast} ((h_{\mbf{x} })_{L',V'}) \frac{}{} \right\} $.
Let
\begin{align}
&H_{L,V,L',V'} \;  = \; \; H_{0,V,V'} \; + \; \kappa H_{L,V,L' , V'}' ,   \\
&H_{L,V,L'} \;\; \; \;  = \; \; H_{0,V} \; + \; \kappa H_{L,V, L'}'  .
\end{align}
In a similar way to  $H$, we can prove that $H_{L,V,L'} $ and $H_{L,V,L' ,V'} $  are
 essentially self-adjoint on $\ms{D}_{0,V}$ and $\ms{D}_{0,V,V'} = \YukawacoreVVprime$, respectively.

\begin{lemma}   \label{11/8.3}

Suppose   \textbf{(A.1)}-\textbf{(A.3)}.  Then   $H_{L,V,L',V' } $ is reduced to
  $\ms{F}_{\Yukawa  V, V' } $, and
 $H_{L,V,L',V' \restr \, \ms{F}_{\Yukawa  V, V' }} $ has  purely discrete spectrum  in $ [ E_{0} ( H_{L, V, L' ,V'}  ), $
$E_{0}( H_{L,V,L',V'}   )+ \nu  ) $.
\end{lemma}
\textbf{(Proof)}
In a similar way to the proof of Lemma \ref{11/8.1}, it is shown that $H_{L,V,L' ,V'} $ is reduced to
 $\ms{F}_{\Yukawa V,V'} $. Since $H_{0,V,V' \restr \ms{F}_{\Yukawa V,V'}} $ has a compact resolvent,
  $H_{L,V,L',V' \restr \, \ms{F}_{\Yukawa  V, V' }} $ also     has a  compact resolvent by
   the general theorem \cite[Theorem 3.8]{AH97}.
   Hence, in particular, $H_{L,V,L',V' \restr \, \ms{F}_{\Yukawa  V, V' }} $ has  purely discrete spectrum  in $ [ E_{0} ( H_{L, V, L' ,V'}  ), $
$E_{0}( H_{L,V,L',V'}   )+ \nu  ) $. $\blacksquare$

\begin{lemma} \label{11/8.4}
Assume \textbf{(A.1)}-\textbf{(A.4)}.
\label{resolventconvergence}
Then for all $ z \in \mbf{C} \setminus \mbf{R} $, it follows that
\[
\mbf{(1)} \; \lim_{V' \to \infty} \| (H_{L,V,L',V' }  -z )^{-1}
- ( H_{L,V,L' } -z )^{-1} \| = 0,
  \qquad \mbf{(2)} \; \lim_{L' \to \infty} \| (H_{L,V,L' }  -z )^{-1}
- ( H_{L,V } -z )^{-1} \| = 0 .
 \]
 \end{lemma}
\textbf{(Proof)}
We see that
\begin{align}
 &( H_{L,V, L',V'} -z )^{-1} -(H_{L,V,L'} -z )^{-1}  \\
&\qquad   =    (H_{L,V,L',V'}  -z )^{-1}  \left\{ \frac{}{} I \otimes
  (H_{\KG} -H_{\KG ,V'}   ) +
 \kappa
(H_{L,V, L'}' - H_{L,V,L',V'}'  ) \right\}  (H_{L,V, L'}-z )^{-1}  .
\label{10/7.3}
\end{align}
Let
  $ \; C_{V' , m }  = \sqrt{3} \left( \frac{\pi}{V'} \right)^{3}( \frac{1}{2m} +1 )$.
  It is shown in
\cite[Lemma 3.1]{AH97} that
\begin{equation}
 \| (I \otimes   (H_{\KG} -H_{\KG ,V'}  ) )  (H_{L,V, L'} -z )^{-1} \|
\leq \frac{2 C_{V',m}}{(1-C_{V' ,m })} \| (I \otimes H_{\KG}) (H_{L,V, L'} -z)^{-1} \| \to 0    \label{10/7.4}
\end{equation}
as $V' \to  \infty $. By (\ref{bounda}) and (\ref{boundad}), we also see that
\begin{align*}
\|   (H_{L,V,L'}' - H_{L,V,L',V'}'  )   (H_{L,V, L'}-z )^{-1} \|
 &\leq \sum_{l,l'} | \gamma^0_{l, l'} | M_{\Dir}^l M_{\Dir}^{l'} \left\{
 \beta_{1} \int_{\Rthree} | \cutoff{I}{x}  |
  \left\|  \frac{(h_{\mbf{x}})_{L'} }{\sqrt{\omega}} -\frac{(h_{\mbf{x}})_{L' , V'} }{\sqrt{\omega}}
\right \| \dx \right.  \\
   & \left. + \beta_{2} \int_{\Rthree}  | \cutoff{I}{x}  | \|  (h_{\mbf{x}})_{L'}  -   (h_{\mbf{x}})_{L', V'} \|  \right\} \dx ,
\end{align*}
where $\beta_{1} =  \sqrt{2}
 \|  I \tens \HKG^{1/2} (H_{L,V,L'} -z )^{-1} \| \; $ and
$\;  \beta_{2}  = \frac{1}{\sqrt{2 }} \|  (H_{L,V,L'} -z )^{-1}  \|$.
From Assumptions  \textbf{(A.2)}, \textbf{(A.4)} and the fact
$| e^{i \mbf{k} \cdot \mbf{x} } - e^{i \mbf{k}' \cdot \mbf{x} } |
 \leq  |\mbf{k} -\mbf{k}' | \,  | \mbf{x}|$,
   it follows that
      $\; \lim\limits_{V' \to \infty}
\int_{\Rthree} |\chi_{\I} (\mbf{x}) | \|  (h_{\mbf{x}})_{L'}  -   (h_{\mbf{x}})_{L', V'} \|  \;  \dx = 0  $
 and $ \; \lim\limits_{V' \to \infty}  \int_{\Rthree} | \cutoff{I}{x}  |
  \left\|  \frac{(h_{\mbf{x}})_{L'} }{\sqrt{\omega}} -\frac{(h_{\mbf{x}})_{L' , V'} }{\sqrt{\omega}}
\right \| \dx  $  $ \; = 0 $.
Hence we have $ \; \lim\limits_{V' \to \infty}
\|   (H_{L,V,L'}' - H_{L,V,L',V'}'  )   (H_{L,V, L'}-z )^{-1}    \|   = 0 $.
Thus we obtain $\mbf{(1)}$. In a similar way to $\mbf{(1)}$,
we can also  prove $\mbf{(2)}$. $\blacksquare$ \\

$\quad$ \\
\textbf{(Proof of Proposition \ref{11/8.2})} \\
The decomposition $
 L^{2} (\Rthree ) =  \ell^{2} ( \Gamma_{V'} )\oplus  \ell^{2} ( \Gamma_{V'} )^{\bot} $ yields that
$\ms{F}_{\KG} \simeq   \ms{F}_{\KG , V'} \otimes (\oplus_{n=0}^{\infty}
 \ell^2 (\Gamma_{V'}  )^{\bot}  )  $. Then we have  $\ms{F}_{\Yukawa  V} \simeq
  \ms{F}_{\Yukawa  V, V' } \oplus (\ms{F}_{\Yukawa  V, V' })^{\bot} $,  where
 $ (\ms{F}_{\Yukawa  V, V' })^{\bot}   =   \oplus_{n=1}^{\infty} \ms{F}_{V, V'}^{(n)}   $ with
  $ \ms{F}_{V, V'}^{(n)} =\ms{F}_{\Yukawa  V, V'} \tens (\otimes_{s}^{n}  \ell^2 (\Gamma_{V'}  )^{\bot}  ) $. Then we have for $n \geq1$, \\
\[
 H_{L,V,L',V' \restr \ms{F}_{V, V'}^{(n)}} \simeq
  H_{L,V,L',V' \; \restr  \ms{F}_{\Yukawa  V, V'}}  \tens I_{ \restr \otimes_{s}^{n}  \ell^2 (\Gamma_{V'}  )^{\bot} }  \; +  I_{\restr  \ms{F}_{\Yukawa V , V'}}  \tens \sqz{\omega}_{\restr \otimes_{s}^{n}  \ell^2 (\Gamma_{V'}  )^{\bot} }  \; \geq E_{0} (H_{L,V, L',V'})  +  n m   .
\]
Hence we have $H_{L,V,L',V' \restr (\ms{F}_{\Yukawa  V, V' })^{\bot}} \geq  E_{0} (H_{L,V, L',V'})  +  \nu $. While
 $H_{L,V,L',V' \restr \, \ms{F}_{\Yukawa  V, V' }} $ has  purely discrete spectrum  in $ [ E_{0} ( H_{L, V, L' ,V'}  ), $
$E_{0}( H_{L,V,L',V'}   )+ \nu  ) $
by Lemma 3.3.
Then
 $H_{L,V,L',V' \restr \, \ms{F}_{\Yukawa  V}} $ also has  purely discrete spectrum  in $ [ E_{0} ( H_{L, V, L' ,V'}  ), $
$E_{0}( H_{L,V,L',V'}   )+ \nu  ) $.
Since $H_{L,V,L',V'}
$ converges to $H_{L,V, L'}$
as $V' \; \to \; \infty$ in the norm resolvent sense
by Lemma \ref{11/8.4},
  $\, H_{L,V, L' }$ has  purely discrete spectrum in  $[ E_{0}(H_{L,V, L'}) ,
E_{0}(H_{L,V, L'} )+ \nu )$ by
\cite[Lemm 4.6]{SH72}.
Since
 $H_{L, V, L'}$ converges to $H_{L,V} $ in the norm resolvent sense
as $L'\;\to \;\infty$ by Lemma \ref{11/8.4},
 $H_{L,V}$ has   also purely discrete spectrum  in  $[ E_{0}(H_{L,V} ) , E_{0} (H_{L,V})+ \nu )$.

\begin{lemma} \label{11/8.5}
Assume \textbf{(A.1)}-\textbf{(A.4)}.
\label{resolventconvergence}
For all $ z \in \mbf{C} \setminus \mbf{R} $, it follows that
\[
\mbf{(1)} \; \lim_{V \to \infty} \| (H_{L,V}  -z )^{-1}
- ( H_{L } -z )^{-1} \| = 0,
 \qquad \mbf{(2)} \; \lim_{L \to \infty} \| (H_{L}  -z )^{-1}
- ( H -z )^{-1} \| = 0 .
 \]
\end{lemma}
\textbf{(Proof) }
The proof is quite parallel with that of Lemma \ref{11/8.4}
 Let
 $ \; C_{V , M}  = \sqrt{3} \left( \frac{\pi}{V} \right)^{3}( \frac{1}{2M +1} )$.
Then
\begin{align}
&\| ( H_{L,V} -z )^{-1} -(H_{L} -z )^{-1}  \| \leq
 \frac{1}{ | \text{Im} z  |} \left\{  \frac{2C_{V ,M  }}{1-C_{V , M } }
  \| (\HDirac \tens I ) (H_{L} -z )^{-1}  \|\right. \notag \\
&   + \left.  C\sum_{l,l'} |\gamma^0_{l,l'} | \int_{\Rthree} |\cutoff{I}{x}|
 (  \frac{}{}
  \| ( \psi_{L,V}^{l} (\mbf{x}) - \psi_{L}^{l} (\mbf{x}) )^{\ast} \| \|  \psi_{L,V }^{l'} (\mbf{x}) \| \; + \;
  \|  \psi_{L}^{l} (\mbf{x}) \| \|  \psi_{L,V}^{l'} (\mbf{x}) - \psi_{L}^{l'} (\mbf{x}) ) \| ) \right\} d \mbf{x} ,
     \end{align}
where $ C  =  \sqrt{2} M_{\KG}^{1} \| (I \tens \HKG^{1/2}) (H_{L} -z)^{-1}  \|
 \; + \; \frac{1}{\sqrt{2}} M_{\KG}^{0} \| (  H_{L} -z)^{-1}  \|  $.
Here we used (\ref{bound_phi}).
 By (\ref{boundbd}) and
   \textbf{(A.4)},  there exists a  constant $c_{L}^l \geq 0  $ such that $ \|  \psi_{L,V }^{l} (\mbf{x}) \| \leq c_L^l$, and
 $ \; \int_{\Rthree} | \chi_{\I} (\mbf{x}) | \|  \psi_{L,V}^{l} (\mbf{x}) - \psi_{L}^{l} (\mbf{x}) ) \|  d \mbf{x} \to 0 $ as $V \to \infty $.
   Then $ \| (H_{L, V}'  -H_{L} ') (H_{L} -z)^{-1} \|  \to \, 0$ as $ V \to \infty $ follows.
 Thus we obtain $\mbf{(1)}$. Similarly we can prove  $\mbf{(2)}$. $\blacksquare$

$\quad$ \\
\textbf{ {\large (Proof  of  Theorem \ref{MainTheorem})}} \\
The proof is parallel with that of Proposition \ref{11/8.2}.
From the decomposition
  $ L^{2} (\Rthree ; \mbf{C}^4 ) =  \ell^{2} ( \Gamma_{V} ; \mbf{C}^4 )
 \oplus  \ell^{2} ( \Gamma_{V} ; \mbf{C}^4 )^{\bot} $,
it follows that
 $\ms{F}_{\Yukawa} \simeq  \FYukawaV \oplus  (\FYukawaV)^{\bot} $,
 where $(\FYukawaV)^{\bot} =\oplus_{n=1}^{\infty}  \ms{F}_{V}^{(n)}$ with   $ \ms{F}_{V}^{(n)} =\FYukawaV \tens (  \natens   \ell^2 (\Gamma_{V} ; \mbf{C}^4)^{\bot}  ) $.
 Then we have for $n \geq1$,
\begin{equation}
 H_{L,V \restr \ms{F}_{V}^{(n)}} \simeq
  H_{L,V \; \restr  \ms{F}_{\Yukawa  V}}  \tens
 I_{ \restr \natens  \ell^2 (\Gamma_{V} ; \mbf{C}^4)^{\bot} }  \; +  I_{\restr  \ms{F}_{\Yukawa V}}  \tens \sqz{\omega}_{\restr \natens \ell^2 (\Gamma_{V} ;   \mbf{C}^4)^{\bot} }  \; \geq E_{0} (H_{L,V})  +  n M   ,  \label{gap_nu}
  \end{equation}
and
   $H_{L,V \restr (\FYukawaV)^{\bot} } $ $\geq  E_{0} (H_{L,V})  +  \nu $.
Then
 $H_{L,V } $ has  purely discrete spectrum  in $ [ E_{0} ( H_{L, V}  ), E_{0}( H_{L,V}   )+ \nu  ) $,
 since
 $H_{L ,V\restr \, \ms{F}_{\Yukawa  V }} $ has  purely discrete spectrum  in $ [ E_{0} ( H_{L, V}  ),
E_{0}( H_{L,V}   )+ \nu  ) $
 by Proposition \ref{11/8.2}.
 Then
Lemma~\ref{11/8.5}
yields that  $H$ has   also purely discrete spectrum  in  $[ E_{0}(H ) , E_{0} (H)+ \nu )$. $\blacksquare $


$\quad $ \\
$\quad $ \\
{\Large \textbf{Acknowledgments}}.  \\
It is pleasure to thank Professor Fumio Hiroshima for his advice and discussions.

\end{document}